\documentclass[a4paper,12pt]{article}
\usepackage{cmap}                        
\usepackage[cp1251]{inputenc}            
\usepackage[english]{babel}
\usepackage[left=2cm,right=2cm,top=2cm,bottom=2cm]{geometry} 
\usepackage{amssymb}
\usepackage{amsmath, amsthm}
\newtheorem{thm}{Theorem}[section]
\newtheorem{cor}[thm]{Corollary}

\newtheorem{prop}[thm]{Proposition}

\newtheorem{defi}[thm]{Definition}

\begin{document}

\begin{center}\Large
\textbf{ARITHMETIC GRAPHS OF  FINITE GROUPS}
\end{center}
\begin{center} V. I. Murashka and A. F. Vasil'ev\end{center}
\begin{center} \{mvimath@yandex.ru, formation56@mail.ru\}\end{center}
\begin{center} Francisk Skorina Gomel State University, Gomel\end{center}

\medskip

\textbf{Abstract.} In this paper we introduced an arithmetic graph function which associates with every group $G$  the directed graph whose  vertices corresponds to the divisors of $|G|$. With the help of such functions we introduced arithmetic graphs of classes of groups, in particular of hereditary saturated formations. We formulated the problem of the recognition of  classes of groups by arithmetic graph functions and investigated this problem  for some arithmetic graph functions.


\textbf{Mathematic Subject Classification}(2010): 20D10, 20D60, 20F17, 05C25.

\medskip



\section{Introduction}

All considered groups are finite. There have been a lot of papers recently in which with every finite group associates certain graph. The considered problem was to analyze the relations between the structure of a group and the properties of its graph. This trend goes back to 1878 when A. Cayley \cite{c1} introduced his graph. This graph has a lot of applications (see \cite{c2}).

Another illustration of this direction is the problem of Paul Erd\"{o}s on non-commuting graphs which was solved in 1976 \cite{c3}.

There exist many different graphs that are connected with a  group. All these graphs can be  divided into three general types. We shall classify them by the type of their vertices.

Graphs of the first type are the graphs whose vertices are elements of a group. Examples of such graphs are the non-commuting graph \cite{c3} and the power graph \cite{f5}.

Graphs of the second type are the graphs whose vertices are subsets of the set of elements of a group. For example the permutability graph \cite{f8} and the intersection graph of subgroups \cite{f9} are such graphs.

Graphs of the third type are the graphs whose vertices are  some  divisors of the order of a group. The Gruenberg-Kegel graph is an example of such graph.

In this paper we will consider only graphs from the third type.

Recall \cite{g3} that the Gruenberg-Kegel or prime graph $\Gamma_p$ of a group $G$ is  the graph with the  set of vertices $V(\Gamma_p(G))=\pi(G)$ where $\pi(G)$ is the set of all prime divisors of $|G|$  and $(p, q)$ is an edge if and only if $G$ contains an element of order $pq$.

This graph is connected to the recognition problem of groups by their graph.
 Recall that a group  $G$ is called recognizable by the prime graph if
 $\Gamma_p(G)=\Gamma_p(H)$ implies $H\simeq G$ for any group $H$.
 There are many non-isomorphic groups with nontrivial solvable radical and the
 same prime graph. Therefore of prime interest
 (for example see \cite{c4}) is this problem  only for simple and almost simple groups. In this paper we will consider the recognition problem up to a class of groups.

\begin{defi} A function $\Gamma:$ $\{$groups$\}\,\rightarrow\,\{$graphs$\}$  is called a graph function.\end{defi}

\begin{defi} Let $\Gamma$ be a graph function and $\mathfrak{X}$ be a class of groups. We shall say that  $\mathfrak{X}$ is  recognized by $\Gamma$ if from $G_1\in\mathfrak{X}$ and $\Gamma(G_1)=\Gamma(G_2)$ it follows that $G_2\in\mathfrak{X}$.      \end{defi}

\textbf{Problem 1.} $(a)$ Let $\Gamma$ be a graph function. Describe all group classes (formations, Fitting classes, Schunk classes) which are recognizable by $\Gamma$.

   $(b)$ Let $\mathfrak{X}$ be a class of groups (formation, Fitting class, Schunk class). Find graph functions $\Gamma$ which recognize $\mathfrak{X}$.

\medskip

In this paper we are interested only in the following graph functions.

\begin{defi}  A graph function $\Gamma$ is called arithmetic  if $V(\Gamma(G))$ is a subset of the set of divisors of  $|G|$  for every group $G$ and  $\Gamma(1)=\emptyset$.  Graph $\Gamma(G)$ is called arithmetic. \end{defi}

The Gruenberg-Kegel graph is the example of an arithmetic graph.

In 1968 T. Hawkes \cite{g5}  considered a directed graph of a group $G$ whose set of vertices is $\pi(G)$ and $(p, q)$ is an edge  if and only if $q\in \pi(G/\mathrm{O}_{p',p}(G))$. We shall call this graph by Hawkes graph and will denote it by $\Gamma_H(G)$. T. Hawkes showed that a group $G$ has a Sylow tower for some linear order $\phi$ if and only if $\Gamma_H(G)$  has no  cycles.
 Note that $\Gamma_H(G)$ might have a loop.



 A. D'Aniello, C. De Vivo and G. Giordano  \cite{g4} introduced the Sylow graph and used it in the studying of properties of lattice formations.  Recall that the Sylow graph $\Gamma_s(G)$ of a group $G$ is the directed graph whose set of vertices is $\pi(G)$ and $(p, q)$ is an edge of $\Gamma_s(G)$ if and only if $q\in \pi(N_G(P)/PC_G(P))$ for a Sylow $p$-subgroup $P$ of $G$.

 Let us note that $q\in \pi(N_G(P)/PC_G(P))$ is equivalent to   $q\in \pi(N_G(P)/\mathrm{O}_{p', p}(N_G(P)))$ for a Sylow $p$-subgroup $P$ of $G$.
From $\mathrm{O}_{p',p}(G)\cap N_G(P)\leq \mathrm{O}_{p',p}(N_G(P))$  it follows that $\Gamma_s(G)\subseteq\Gamma_H(G)$  for any group $G$.

We shall consider new example of an arithmetic graph of a group. Recall that a Schmidt $(p, q)$-group is a Schmidt group $G$ with $\pi(G)=\{p, q\}$ and  the normal Sylow $p$-subgroup.

 \begin{defi} We shall call $\Gamma_{Sch}(G)$ the Schmidt graph of a group $G$ if $V(\Gamma_{Sch}(G))=\pi(G)$  and $(p,q)$ is an edge of $\Gamma_{Sch}(G)$ if and only if $G$ contains a Schmidt $(p, q)$-group as subgroup.  \end{defi}

Recall that a group of automorphisms of a subgroup  $H$ induced by  a group $G$ is $N_G(H)/C_G(H)$.

Note that $(p, q)$  is an edge of $\Gamma_H(G)$ if and only if some element of $G$ acts as nontrivial automorphism of order $q^\alpha$ on some chief  factor $H/K$ where $p\in\pi(H/K)$.

It is easy to see that $(p, q)$  is an edge of $\Gamma_s(G)$ if and only if some element of $G$ acts as nontrivial automorphism of order $q^\alpha$ on a Sylow $p$-subgroup $P$ of $G$.

  Let $\Gamma$ be an arithmetic graph function such that for every group $G$ $V(\Gamma(G))=\pi(G)$ and $(p, q)\in E(\Gamma(G))$ if and only if  $q\in\pi(N_G(P)/PC_G(P))$  for some $p$-subgroup $P$ of $G$. Then $\Gamma=\Gamma_{Sch}$. It is clear that $\Gamma_{Sch}(G)\subseteq\Gamma(G)$ for every group $G$.  Let us show the inverse inclusion. If $(p, q)\in E(\Gamma(G))$ then there is a $p$-subgroup $P$ of $G$  such that $q\in\pi(N_G(P)/C_G(P))$. Let $Q$ be a Sylow $q$-subgroup of $N_G(P)$. Note that $PQ$ is a non-nilpotent biprimary group  with the normal Sylow $p$-subgroup. Hence it contains a Schmidt $(p, q)$-group. So  $(p, q)\in E(\Gamma_{Sch}(G))$. Thus  $\Gamma=\Gamma_{Sch}$.

This lead us to the following   construction of arithmetic graph functions.
Recall that a section of a group $G$ is $H/K$ where $H\leq G$ and $K\triangleleft H$.

\begin{defi} We shall call a function $f$ the section functor if $f$   assigns to each group $G$ a possibly empty set  $f(G)$ of sections of $G$ satisfying $\alpha(f(G))=f(\alpha(G))$ for any isomorphism $\alpha: G\rightarrow G^*$. \end{defi}

\begin{defi} Let $\theta$ be a function which assigns to every prime $p$ a section functor $\theta(p)$.

  $(1)$ We shall call $\Gamma_\theta$  the $\theta$-local arithmetic graph function if   $V(\Gamma_\theta(G))=\pi(G)$ and $(p, q)\in E(\Gamma_\theta(G))$  if and only if  $q\in\pi(\{ N_G(P)/C_G(P) | P\in\theta(p)(G)\})$ for every group $G$.

   $(2)$ We shall call $\Gamma_{\overline{\theta}}$  the $\overline{\theta}$-local arithmetic graph function if   $V(\Gamma_{\overline{\theta}}(G))=\pi(G)$ and $(p, q)\in E(\Gamma_{\overline{\theta}}(G))$  if and only if $q\neq p$ and $q\in\pi(\{ N_G(P)/C_G(P) | P\in\theta(p)(G)\})$ for every group $G$.
  \end{defi}

 If $\theta(p)$ is the set of all chief factors $H/K$  of $G$ such that $p\in\pi(H/K)$ then $\Gamma_{\theta}=\Gamma_H$.

 If $\theta(p)$ is the set of all Sylow subgroups of $G$ then $\Gamma_{\overline{\theta}}=\Gamma_s$.

 If $\theta(p)$ is the set of all $p$-subgroups of $G$ then $\Gamma_{\overline{\theta}}=\Gamma_{Sch}$.

It is straightforward to check that

\begin{prop} $\Gamma_s(G)\subseteq\Gamma_{Sch}(G)\subseteq\Gamma_H(G)$ for every group $G$.  \end{prop}

  \section{Preliminary results}

 We use standard terminology and notation on groups and graphs that can be found in \cite{s7, s8, gr}. Recall that $\pi(\mathfrak{F})=\underset{G\in\mathfrak{F}}\cup\pi(G);$. A Schmidt group is a non-nilpotent group all whose proper subgroups are nilpotent.

Recall that for a class of group $\mathfrak{X}$
\[\textbf{Q}\mathfrak{X}=(G| \exists H\in\mathfrak{X} \, and \, an\, epimorphism \, from \, H \, onto \,\, G)\]
\[\textbf{R}_0\mathfrak{X}=(G| \exists N_i\triangleleft G  \,\, (i=1,\dots, n) \,with\, G/N_i\in\mathfrak{X} \,  and \, \cap_{i=1}^nN_i=1)\]
\[\textbf{S}\mathfrak{X}=(G| \exists H\in\mathfrak{X} \, and \,  G\leq H)\]
\[\textbf{D}_0\mathfrak{X}=(G| \exists N_i\in\mathfrak{X}   \,\, (i=1,\dots, n) \,with\, G=N_1\times\dots\times N_n)\]
\[\textbf{N}_0\mathfrak{X}=(G| \exists N_i\in\mathfrak{X} \,\, and \,\, N_i\triangleleft G  \,\,  (i=1,\dots, n) \,with\, G= N_1\cdot\cdot\cdot N_n)\]
\[\textbf{E}_\Phi\mathfrak{X}=(G| \exists H\in\mathfrak{X}   \,with\, G/\Phi(G)\simeq H)\]
\[\textbf{DM}\mathfrak{X}=(G| \exists K\in\mathfrak{X}   \,\,   \,with\, K=G\times H)\]

A class $\mathfrak{F}=\{\textbf{R}_0, \textbf{Q}\}\mathfrak{F}$   is called a formation. An $\textbf{S}$-closed formation is called hereditary. It is well known that if $\mathfrak{F}=\textbf{QSD}_0\mathfrak{F}$ then $\mathfrak{F}$ is a hereditary formation. An $\textbf{E}_{\Phi}$-closed formation is called saturated.

A function $f$: $\mathbb{P}\rightarrow\{formations\}$ is called a formation function.   Formation  $\mathfrak{F}=LF(f)$ defined by a formation function $f$: $LF(f)=(
G|$ if $H/K$ is a chief factor of $G$ and $p\in\pi(H/K)$ then $G/C_G(H/K)\in
f(p))$ is called local. By well known  Gashutz-Lubeseder-Schmid Theorem local formations are exactly saturated formations.

A formation $\mathfrak{F}$ is called formation with the Shemetkov property is every $s$-critical for $\mathfrak{F}$ group is either a Schmidt group or a group of prime order.

A formation $\mathfrak{F}$ is called solubly saturated if from $G/\Phi(G_\mathfrak{S})\in\mathfrak{F}$ it follows that $G\in\mathfrak{F}$ where $G_\mathfrak{S}$ is the soluble radical of $G.$

Recall \cite{gr} that here a graph $\Gamma$ is a pair of sets $V$ and $E$ where $V$ is a set of vertices of $\Gamma$ and $E$ is a set of edges of $\Gamma$, i.e. the set of ordered pairs of elements from $V$. An edge $(v, v)$ is called a loop. If $(a, b)\in E(\Gamma)$ if and only if $(b, a)\in E(\Gamma)$ for any $(a, b)\in E(\Gamma)$ then a graph $\Gamma$ is called undirected. For a given graph $\Gamma$ $V(\Gamma)$ is the set of vertices of $\Gamma$ and   $E(\Gamma)$ is the set of edges of $\Gamma$. Two graphs $\Gamma_1$ and $\Gamma_2$  are called equal ($\Gamma_1=\Gamma_2$) if $V(\Gamma_1)=V(\Gamma_2)$ and $E(\Gamma_1)=E(\Gamma_2)$. Graph $\Gamma_1$ is called subgraph of $\Gamma_2$   ($\Gamma_1\subseteq\Gamma_2$) if $V(\Gamma_1)\subseteq V(\Gamma_2)$ and $E(\Gamma_1)\subseteq E(\Gamma_2)$. Graph $\Gamma$ is called a union of graphs  $\Gamma_1$ and $\Gamma_2$ ($\Gamma=\Gamma_1\cup\Gamma_2$) if $V(\Gamma)=V(\Gamma_1)\cup V(\Gamma_2)$ and $E(\Gamma)=E(\Gamma_1)\cup E(\Gamma_2)$.



\section{Arithmetic graphs and classes of groups}

The notion of arithmetic graph function is very wide. That is why  we are interested only in those arithmetic graph functions  which are connected with closure operators on classes of groups.

\begin{defi} An arithmetic graph function $\Gamma$ is called:

$(a)$ $\textbf{S}$-closed  if $H\leq G$ implies $\Gamma(H)\subseteq\Gamma(G)$.

$(b)$  $\textbf{Q}$-closed  if  $\Gamma(G/N)\subseteq\Gamma(G)$ for any group $G$ and normal subgroup $N$ of  $G$.

$(c)$ $\textbf{D}_0$-closed  if $\Gamma(G_1\times\dots\times G_n)=\Gamma(G_1)\cup\dots\cup\Gamma(G_n)$ for any groups $G_1,\dots, G_n$.

$(d)$ $\textbf{R}_0$-closed if  $\Gamma(G/N_1)\cup\dots\cup\Gamma(G/N_n)=\Gamma(G/\cap_{i=1}^nN_i)$ for any group $G$ and its normal subgroups $N_1,\dots, N_n$.

$(e)$  $\textbf{N}_0$-closed if  $\Gamma(N_1)\cup\dots\cup\Gamma(N_n)=\Gamma(G)$ for any group $G=N_1\cdot\cdot\cdot N_n$ and its normal subgroups $N_1,\dots, N_n$.

$(f)$ $\textbf{E}_\Phi$-closed if $\Gamma(G/\Phi(G))=\Gamma(G)$ for any group $G.$

\end{defi}

\emph{Remark.} If we want to verify that an arithmetic graph function is $\textbf{C}$-closed where $\textbf{C}\in\{\textbf{D}_0, \textbf{R}_0, \textbf{N}_0\}$ we need only to consider the case when $n=2$.

Let us note that Gruenberg-Kegel graph function is $\textbf{S}$-closed and $\textbf{Q}$-closed but not  $\textbf{D}_0$-closed.

\begin{prop}\label{p1} $\Gamma_H$ is  $\{\textbf{S}, \textbf{Q}, \textbf{D}_0, \textbf{R}_0, \textbf{N}_0, \textbf{E}_\Phi\}$-closed graph function.\end{prop}

\emph{Proof.} Let $\Gamma=\Gamma_{H}$. Let us show that $\Gamma$ is $\textbf{S}$-closed. Let $H$ be a subgroup of a group $G$.  From $H\cap \mathrm{O}_{p', p}(G)\leq \mathrm{O}_{p', p}(H)$ it follows that $H/\mathrm{O}_{p', p}(H)$ is  the quotient group of $H/H\cap \mathrm{O}_{p', p}(G)\simeq H\mathrm{O}_{p', p}(G)/\mathrm{O}_{p', p}(G)$. So    $\Gamma(H)\subseteq\Gamma(G)$. Hence $\Gamma$ is $\textbf{S}$-closed.

Let us show that $\Gamma$ is $\textbf{Q}$-closed. From $\mathrm{O}_{p', p}(G)N/N\leq \mathrm{O}_{p', p}(G/N)$ it follows that $G/N/\mathrm{O}_{p', p}(G/N)$ is the quotient group of $G/N/\mathrm{O}_{p', p}(G)N/N\simeq G/\mathrm{O}_{p', p}(G)N$. Thus    $\Gamma(G/N)\subseteq\Gamma(G)$. Hence $\Gamma$ is $\textbf{Q}$-closed.

Let us show that $\Gamma$ is $\textbf{D}_0$-closed. Let $G_1$ and $G_2$ be a groups. Note that $G_1\leq G_1\times G_2$. Since $\Gamma$ is $\textbf{S}$-closed, we see that $\Gamma(G_1)\subseteq \Gamma(G_1\times G_2)$. By analogy $\Gamma(G_2)\subseteq \Gamma(G_1\times G_2)$. Hence $\Gamma(G_1)\cup\Gamma(G_1)\subseteq \Gamma(G_1\times G_2)$.

Note that $\mathrm{O}_{p', p}(G_1)\times \mathrm{O}_{p', p}(G_1)\leq \mathrm{O}_{p', p}(G_1\times G_2)$. Now $G_1\times  G_2/\mathrm{O}_{p', p}(G_1\times G_2)$ is the quotient group of $G_1\times G_2/\mathrm{O}_{p', p}(G_1)\times \mathrm{O}_{p', p}(G_1)\simeq G_1/\mathrm{O}_{p', p}(G_1)\times G_2/\mathrm{O}_{p', p}(G_2)$. So $\Gamma(G_1)\cup\Gamma(G_1)\supseteq \Gamma(G_1\times G_2)$. Now $\Gamma(G_1)\cup\Gamma(G_1)= \Gamma(G_1\times G_2)$. Thus $\Gamma$ is $\textbf{D}_0$-closed.

Let us show that $\Gamma$ is $\textbf{R}_0$-closed. Let  $A$ and $B$ be normal subgroups of a group $G$.    It is clear that  $\Gamma(G)\supseteq\Gamma(G/A)\cup\Gamma(G/B)$. Note that $G$ isomorphic to a subgroup of $G/A\times G/B$. Hence $\Gamma(G)\subseteq\Gamma(G/A)\cup\Gamma(G/B)$. Thus $\Gamma(G)=\Gamma(G/A)\cup\Gamma(G/B)$.

Let us show that $\Gamma$ is $\textbf{N}_0$-closed. Let $G=AB$ where $A$ and $B$ are normal subgroups of $G$.    It is clear that  $\Gamma(G)\supseteq\Gamma(A)\cup\Gamma(B)$. Note that $\mathrm{O}_{p', p}(A)\mathrm{O}_{p', p}(B)\leq \mathrm{O}_{p', p}(G)$, $\mathrm{O}_{p', p}(A)\cap \mathrm{O}_{p', p}(A)\mathrm{O}_{p', p}(B)=\mathrm{O}_{p', p}(A)$ and $\mathrm{O}_{p', p}(A)\mathrm{O}_{p', p}(B)\cap \mathrm{O}_{p', p}(B)=\mathrm{O}_{p', p}(B)$. Now $G/\mathrm{O}_{p', p}(G)$ is a quotient group of $G/\mathrm{O}_{p', p}(A)\mathrm{O}_{p', p}(B)\simeq (A/\mathrm{O}_{p', p}(A)) (B/\mathrm{O}_{p', p}(B))$. Hence   $\Gamma(G)\subseteq\Gamma(A)\cup\Gamma(B)$. Thus $\Gamma(G)=\Gamma(A)\cup\Gamma(B)$.

It is well known that $\Phi(G)\leq \mathrm{O}_{p', p}(G)$ and $\mathrm{O}_{p', p}(G/\Phi(G))=\mathrm{O}_{p', p}(G)/\Phi(G)$. Hence $\Gamma$ is $\textbf{E}_{\Phi}$-closed. $\square$

\begin{prop}\label{p2} $\Gamma_{Sch}$ is  $\{\textbf{S}, \textbf{Q}, \textbf{D}_0, \textbf{R}_0\}$-closed graph function and is not $\textbf{E}_\Phi$-closed.\end{prop}

\emph{Proof.} Let $\Gamma=\Gamma_{Sch}$. It is evident that $\Gamma$ is   $\textbf{S}$-closed. Let $G$ be a group such that it does not have Schmidt $(p, q)$-subgroup. It means that all $\{p, q\}$-subgroups of $G$ are $q$-nilpotent. So for every normal subgroup $N$ of $G$  all $\{p, q\}$-subgroups of $G/N$ are $q$-nilpotent. Hence     $G/N$ does not have Schmidt $(p, q)$-subgroup. It means that $\Gamma_{Sch}(G/N)\subseteq\Gamma_{Sch}(G)$. So  $\Gamma$ is   $\textbf{Q}$-closed.

Let $G_1$ and $G_2$ be a groups. Assume that $\Gamma_{Sch}(G_1)\cup\Gamma_{Sch}(G_2)\subset\Gamma_{Sch}(G_1\times G_2)$. It means that $G_1$ and $G_2$ does not contain any Schmidt $(p, q)$-groups   and $G_1\times G_2$ contains a Schmidt $(p, q)$-group $H$. Let $\rho_i$ be the projection $H$ on $G_i$ for $i=1, 2$. Then $ker\rho_1$ contains the Sylow $p$-subgroup of $H$. By analogy  $ker\rho_2$ contains the Sylow $p$-subgroup of $H$. Therefore $ker\rho_1\cap ker\rho_2\neq1$, a contradiction. Thus $\Gamma_{Sch}$ is $\textbf{D} _0$-closed.

By analogy with the proof of proposition \ref{p1} it is easy to show that $\Gamma$ is $\textbf{R}_0$-closed.

Let $G\simeq A_5$  the alternating group  of degree 5. Let us note that $G$  does not contain Schmidt $(3, 5)$-subgroup. According to \cite{20} there is faithful
 irreducible Frattini $\mathbb{F}_3G$-module $A$. By known Gaschutz
 theorem \cite{41}, there
exists a Frattini extension  $A\rightarrowtail R\twoheadrightarrow G$
such that $A\stackrel {G}{\simeq} \Phi(R)$ and $R/\Phi(R)\simeq G$. Let $K$ be a subgroup of $R$ such that $K/\Phi(R)\simeq Z_5$ the cyclic group of order 5. Since $A$ is faithful
   $\mathbb{F}_3G$-module, we see that $K$ is not nilpotent group with normal Sylow $3$-subgroup. Hence it contains a Schmidt $(3, 5)$-subgroup. Now $\Gamma(R/\Phi(R))\subset\Gamma(R)$.          $\square$

\smallskip

\begin{prop}\label{p3} $\Gamma_s$ is  $\{\textbf{Q},  \textbf{R}_0\}$-closed graph function and is not $\textbf{S}$-closed.\end{prop}

\emph{Proof.} Let $\Gamma=\Gamma_s$. It is well known that $N_G(P)N/N=N_{G/N}(PN/N)$ for every Sylow subgroup $P$ and normal subgroup $N$ of a group $G$. Let $aN$ be a $q$-element from $N_{G/N}(PN/N)\setminus C_{G/N}(PN/N)$. It means that there is $b\in P$ such that $[a, b]\not\in N $.  Without loose of generality we may assume that $a$ is $q$-element from $N_G(P)$. Note that $[a, b]\neq 1$. So $a\in N_G(P)\setminus C_G(P)$. Hence $\Gamma(G/N)\subseteq \Gamma(G)$. Thus  $\Gamma$ is $\textbf{Q}$-closed.

Let us show that $\Gamma$ is $\textbf{R}_0$-closed. Let $N_1$ and $N_2$ be a normal subgroups of   $G$ such that $N_1\cap N_2=1$. It is clear that $\Gamma(G/N_1)\cup \Gamma(G/N_2)\subseteq \Gamma(G)$. Let $P$ be a Sylow $p$-subgroup of $G$ and   $a\in N_G(P)\setminus C_G(P)$ be a $q$-element of $G$. It means that $[a, b]\neq 1$ for some $b\in P$. Assume that $aN_i\in C_{G/N_i}(PN_i/N_i)$ for $i=1,2$. Hence $[a, b]\in N_1\cap N_2=1$, a contradiction. So  $\Gamma(G/N_1)\cup \Gamma(G/N_2)\supseteq \Gamma(G)$. Thus $\Gamma(G/N_1)\cup \Gamma(G/N_2)= \Gamma(G)$ and $\Gamma$ is $\textbf{R}_0$-closed.

The Sylow graph of the symmetric group of degree 4 is $(3, 2)$. The Sylow graph of the alternating group of degree 4 is $(2, 3)$. So $\Gamma$ is not $\textbf{S}$-closed. $\square$

\smallskip

\textbf{Problem 2.} Is $\Gamma_{Sch}$ $\textbf{N}_0$-closed? Is $\Gamma_{s}$ $\textbf{N}_0$-closed?  Is $\Gamma_{s}$ $\textbf{E}_{\Phi}$-closed?

\medskip

Now we will introduce a graph of a class of groups.

\begin{defi} Let $\mathfrak{X}$ be a class of groups and $\Gamma$   be an arithmetic graph function. Then $\Gamma(\mathfrak{X})=\bigcup\limits_{G\in\mathfrak{X}}\Gamma(G)$.\end{defi}

Let us mention that\[\Gamma(\mathfrak{X}\cup\mathfrak{H})=\bigcup\limits_{G\in\mathfrak{X}\cup\mathfrak{H}}\Gamma(G)=
\bigcup\limits_{G\in\mathfrak{X}}\Gamma(G)\cup\bigcup\limits_{G\in\mathfrak{H}}\Gamma(G)=
\Gamma(\mathfrak{X})\cup\Gamma(\mathfrak{H})\]

 Let us note that $\Gamma(G)$ has a finite set of vertices but $\Gamma(\mathfrak{X})$ may have an infinite set of vertices.   Every subgraph of $\Gamma(\mathfrak{X})$ with a finite set of vertices lies in the finite join of $\Gamma(G_i)$ for some groups $G_i\in\mathfrak{X}$.

\medskip

The arithmetic graph function $\Gamma$ defines the equivalence relation $\backsim$ on the set $\mathcal{G}$ of all classes of groups in the following way. We shall call classes of groups $\mathfrak{X}$ and $\mathfrak{H}$ equivalent  if $\Gamma(\mathfrak{X})=\Gamma(\mathfrak{H})$. If $\mathfrak{X}\backsim\mathfrak{H}$ then $\mathfrak{X}\backsim\mathfrak{H}\cup\mathfrak{X}$. So from Zorn's lemma it follows that every equivalence class has the greatest element by inclusion. The constructive description of this element is given by


\begin{prop} Let $\mathcal{K}\in \mathcal{G}/\backsim$ and $\mathfrak{X}\in\mathcal{K}$. Then $\mathfrak{X}_\Gamma=(G| \Gamma(G)\subseteq \Gamma(\mathfrak{X}))$ is the greatest element of $\mathcal{K}$.   \end{prop}

\emph{Proof.} Let $\mathfrak{H}\in\mathcal{K}$ and $\mathfrak{H}\not\subseteq\mathfrak{X}_\Gamma$. Then there is a group $G\in \mathfrak{H}\setminus\mathfrak{X}_\Gamma$ such that $\Gamma(G)\subseteq\Gamma(\mathfrak{H})=\Gamma(\mathfrak{X})$. Now $G\in\mathfrak{X}_\Gamma$, the contradiction. $\square$


\begin{thm}\label{t1} Let $\textbf{C}\in\{\textbf{S}, \textbf{Q}, \textbf{D}_0, \textbf{R}_0, \textbf{N}_0, \textbf{E}_\Phi\}$ and $\Gamma$ is an arithmetic graph function. If $\Gamma$ is $\textbf{C}$-closed then $\mathfrak{X}_\Gamma$ is $\textbf{C}$-closed for any class of groups $\mathfrak{X}$.\end{thm}

\emph{Proof.} Let $\textbf{C}\in\{\textbf{S}, \textbf{Q}, \textbf{D}_0, \textbf{R}_0, \textbf{N}_0, \textbf{E}_\Phi\}$ and $\Gamma$ is an arithmetic graph function.  From the definition of $\textbf{C}$ for classes of groups and graph functions we see that the graphs of groups from $\textbf{C}\mathfrak{X}$ are either unions of graphs of groups from $\mathfrak{X}$ or subgraphs of graphs of groups from $\mathfrak{X}$. Hence  $\bigcup\limits_{G\in\mathfrak{X}}\Gamma(G)=\bigcup\limits_{G\in\textbf{C}\mathfrak{X}}\Gamma(G)$.
Thus $\Gamma(\mathfrak{X})=\Gamma(\textbf{C}\mathfrak{X})$. It means that $\textbf{C}\mathfrak{X}\subseteq\mathfrak{X}_\Gamma$. Hence $\textbf{C}\mathfrak{X}_\Gamma\subseteq\mathfrak{X}_\Gamma$. Thus $\textbf{C}\mathfrak{X}_\Gamma=\mathfrak{X}_\Gamma$. $\square$



The following proposition gives another definition of $\Gamma_H(\mathfrak{F})$ for a local formation $\mathfrak{F}$.

 \begin{prop}\label{38} Let $f$ be the minimal local definition of local formation $\mathfrak{F}$. Let $\Gamma(\mathfrak{F})$ be a directed graph with $V(\Gamma(\mathfrak{F}))=\pi(\mathfrak{F})$ and $(p,q)\in E(\Gamma(\mathfrak{F}))$ if and only if $q\in\pi(f(p))$. Then $\Gamma(\mathfrak{F})=\Gamma_H(\mathfrak{F})$.\end{prop}

\emph{Proof.} Let us note that $V(\Gamma(\mathfrak{F}))=V(\Gamma_H(\mathfrak{F}))$. From $q\in\pi(f(p))$ it follows that there exists $G\in\mathfrak{F}$ with $q\in\pi(G/\mathrm{O}_{p', p}(G))$. Hence    $\Gamma(\mathfrak{F})\subseteq\Gamma_H(\mathfrak{F})$.

Assume that $(p, q)\in E(\Gamma_H(\mathfrak{F}))$. So there exists $G\in\mathfrak{F}$ with $q\in\pi(G/\mathrm{O}_{p', p}(G))$. From $ G/\mathrm{O}_{p', p}(G)\in f(p)$ it follows that $q\in\pi(f(p))$. Hence    $\Gamma(\mathfrak{F})\supseteq\Gamma_H(\mathfrak{F})$. Thus   $\Gamma(\mathfrak{F})=\Gamma_H(\mathfrak{F})$. $\square$

\smallskip

Let $\pi$ be a set of primes. Recall that $\mathfrak{G}_\pi$ is the class of all $\pi$-groups. If $\pi$ is empty then we assume that $\mathfrak{G}_\pi=(1)$ is the class of unit groups.

\begin{prop}\label{3.9} Let $\mathfrak{F}$ be a class of groups and $\sigma(p)=\{ q | (p, q)\in E(\Gamma_H(\mathfrak{F}))\}$. Then $\mathfrak{F}_{\Gamma_H}=LF(f)$ where  $f(p)=\mathfrak{G}_{\sigma(p)}$ for $p\in\pi(\mathfrak{F})$ and $f(p)=\emptyset$ otherwise.\end{prop}

\emph{Proof.} Let $\mathfrak{F}^*=LF(f)$ where  $f(p)=\mathfrak{G}_{\sigma(p)}$ for $p\in\pi(\mathfrak{F})$ and $f(p)=\emptyset$ otherwise. By proposition \ref{38} $\mathfrak{F}^*\subseteq\mathfrak{F}_{\Gamma_H}$. Assume that $\mathfrak{F}^*\subset\mathfrak{F}_{\Gamma_H}$. Let $G\in\mathfrak{F}_{\Gamma_H}\setminus\mathfrak{F}^*$. It means that there is $p\in\pi(G)$ such that $\pi(G/\mathrm{O}_{p', p}(G))\not\subseteq\sigma(p)$. Hence $\Gamma_H(G)\not\subseteq\Gamma_H(\mathfrak{F})$, a contradiction. Thus $\mathfrak{F}^*=\mathfrak{F}_{\Gamma_H}$. $\square$

\begin{prop}\label{310} Let $\mathfrak{F}$ be a class of groups.  Then $\mathfrak{F}_{\Gamma_{Sch}}$ is a hereditary solubly saturated formation with the Shemtkov property.   \end{prop}

\emph{Proof.} Let $\Gamma=\Gamma_{Sch}$. From proposition \ref{p1} and theorem \ref{t1} it follows that $\mathfrak{F}_\Gamma$ is hereditary formation. Let $G$ be a $s$-critical  group for $\mathfrak{F}_\Gamma$. If $G$ is not a Schmidt group then $E(\Gamma(G))\subseteq E(\Gamma(\mathfrak{F}))$. It means that  $V(\Gamma(G))\not\subseteq V(\Gamma(\mathfrak{F}))$. Hence $G\simeq Z_q$ for some $q\in\mathbb{P}$. Thus $\mathfrak{F}$ is a formation with the Shemetkov property. According to corollary 2.4.3 \cite[p. 94]{s7}   every hereditary formation with the Shmetkov property is solubly saturated.

Assume that $\mathfrak{F}$ is hereditary composition formation with the Shemetkov property and $\mathfrak{F}\neq\mathfrak{F}_\Gamma$, i.e  $\mathfrak{F}\subseteq\mathfrak{F}_\Gamma$. Let $G$ be a group of minimal order from $\mathfrak{F}_\Gamma\setminus\mathfrak{F}$. So $G$ is $s$-critical for $\mathfrak{F}$ and is a Schmidt group.    Since $\Gamma(\mathfrak{F}_\Gamma)=\Gamma(\mathfrak{F})$, we see that $\Phi(G)\neq 1$ and $G/\Phi(G)\in\mathfrak{F}$. Note that $G/\Phi(G)\in\mathfrak{F}\cap\mathfrak{S}$ and  $\mathfrak{F}\cap\mathfrak{S}$  is saturated. Hence $G\in\mathfrak{F}\cap\mathfrak{S}\subseteq\mathfrak{F}$, a contradiction. Thus $\mathfrak{F}=\mathfrak{F}_\Gamma$. $\square$

\begin{prop}Let $\mathfrak{F}$ be a class of groups. Then $\mathfrak{F}_{\Gamma_s}$ is a formation.  \end{prop}

\emph{Proof.} Follows from theorem \ref{t1} and proposition \ref{p3}.$\square$

\section{Classes of groups that are recognized by graphs}

\begin{thm} Let $\mathfrak{F}$ be a hereditary formation and $\pi=\pi(\mathfrak{F})$. Then $\mathfrak{F}$ is  recognized by $\Gamma_p$ if and only if $\mathfrak{F}=\mathfrak{G}_\pi$.  \end{thm}

\emph{Proof.} Let $\mathfrak{F}$ be  recognized by $\Gamma_p$. Since $\mathfrak{F}$ is a hereditary, we see that $Z_p\in\mathfrak{F}$ for all $p\in\pi$. Let $G=\underset{p\in\pi}\times Z_p$. Then $G\in\mathfrak{F}$ and $\Gamma_p(G)$ is the full graph on $\pi$. Let $H$ be a $\pi$-group. Then $\Gamma_p(G)=\Gamma_p(G\times H)$. So $G\times H\in\mathfrak{F}$. Hence $H\in\mathfrak{F}$. Thus $\mathfrak{F}=\mathfrak{G}_\pi$.

If $\mathfrak{F}=\mathfrak{G}_\pi$. It is clear that $\mathfrak{F}$ contains  every group $G$ with $V(\Gamma_p(G))\subseteq\pi$. So $\mathfrak{F}$ is  recognized by $\Gamma_p$. $\square$

\medskip

The following theorem is the main result of this paper.

\begin{thm}\label{4.2} Let $\Gamma$ be a  $\textbf{D}_0$-closed  arithmetic graph function and $\mathfrak{X}$ be a  $\{\textbf{D}_0, \textbf{DM}\}$-closed class of groups. Then $\mathfrak{X}$ is recognized by $\Gamma$  if and only if $\mathfrak{X}=\mathfrak{X}_\Gamma$.
\end{thm}

\emph{Proof.} Assume that $\Gamma$ recognizes $\mathfrak{X}$ and $\mathfrak{X}\neq\mathfrak{X}_\Gamma$. Since $\Gamma(\mathfrak{X})=\Gamma(\mathfrak{X}_\Gamma)$, there exist $\mathfrak{X}$-groups $G_1, \dots, G_n$ such that $\Gamma(G)\subseteq\bigcup\limits_{i=1}^n\Gamma(G_i)$ for every $G\in\mathfrak{X}_\Gamma$. Now $G_1\times\dots\times G_n\in \mathfrak{X}$ and $\Gamma(G_1\times\dots\times G_n)=\Gamma(G\times G_1\times\dots\times G_n)$. Hence $G\times G_1\times\dots\times G_n\in\mathfrak{X}$. So $G\in\mathfrak{X}$. Hence $\mathfrak{X}_\Gamma\subseteq\mathfrak{X}$. From the definition of $\mathfrak{X}_\Gamma$ it follows that $\mathfrak{X}\subseteq\mathfrak{X}_\Gamma$. Thus $\mathfrak{X}=\mathfrak{X}_\Gamma$, a contradiction.

  Since $\mathfrak{X}_\Gamma$ is the maximal class of groups with $\Gamma(\mathfrak{X}_\Gamma)=\Gamma(\mathfrak{X})$, we see that    $\mathfrak{X}_\Gamma$ contains every group $G$ with $\Gamma(G)\subseteq\Gamma(\mathfrak{X}_\Gamma)$. Hence $\Gamma$ recognizes $\mathfrak{X}_\Gamma$. $\square$

\begin{thm}  Let $\mathfrak{F}$ be a  formation and $\sigma(p)=\{ q | (p, q)\in E(\Gamma_H(\mathfrak{F}))\}$. Then $\mathfrak{F}$ is recognized by $\Gamma_H$ if and only if $\mathfrak{F}=LF(f)$  where  $f(p)=\mathfrak{G}_{\sigma(p)}$ for $p\in\pi(\mathfrak{F})$ and $f(p)=\emptyset$ otherwise.\end{thm}

\emph{Proof.} The result follows from theorem \ref{4.2} and proposition \ref{3.9}. $\square$

\begin{thm}  Let $\mathfrak{F}$ be a  formation. Then $\mathfrak{F}$ is recognized by $\Gamma_{Sch}$ if and only $\mathfrak{F}$ is a hereditary solubly saturated formation with the Shemetkov property.\end{thm}

\emph{Proof.} The result follows from theorem \ref{4.2}, proposition \ref{310} and   corollary 2.4.3 \cite[p. 94]{s7}. $\square$

\smallskip

According to proposition \ref{p3} if a formation $\mathfrak{F}$ is recognized by $\Gamma_s$ then it is not hereditary in the general. From theorem \ref{4.2} it follows that   if a hereditary formation $\mathfrak{F}$ is recognized by $\Gamma_s$ then $\mathfrak{F}_{\Gamma_s}$ is hereditary. The following theorem gives additional information about $\mathfrak{F}_{\Gamma_s}$ in this case.

\begin{thm} Let $\mathfrak{F}$ be a hereditary formation. If  $\mathfrak{F}$ is recognized by $\Gamma_s$    then the following statements holds:

$(a)$ $\Gamma_{s}(\mathfrak{F})$ is undirected.

$(b)$  $\mathfrak{F}$ is a solubly saturated formation with the Schemetkov property. \end{thm}

\emph{Proof.} $(a)$ Let $\Gamma=\Gamma_s$ and $(p, q)$ be an edge of $\Gamma(\mathfrak{F})$. Then there is a group $G\in\mathfrak{F}$ such that $(p, q)$ is an edge of $\Gamma(G)$. According to theorem 10.3B \cite{s8} for a cyclic group $Z_q$ there exists a faithful irreducible  $Z_q$-module $T_p$ over the field $\mathbb{F}_p$. Let $H(p, q)=T_p\leftthreetimes Z_q$. Note that $\Gamma_s(H(p, q))$ is $(p, q)$. From $\Gamma_s(G)=\Gamma_s(G\times H(p, q))$ it follows that $G\times H(p, q)\in\mathfrak{F}$. Hence $H(p, q)\in\mathfrak{F}$.

By theorem 10.3B \cite{s8} for a group $H(p, q)$ there exists a faithful irreducible  $H(p, q)$-module $T_q$ over the field $\mathbb{F}_q$. Let $G(p, q)=T_q\leftthreetimes H(p, q)$. Note that every Sylow $q$-subgroup of $G(p, q)$ is a maximal abnormal subgroup of $G(p, q)$. Hence $\Gamma(G(p, q))=\Gamma(H(p, q))$. So $G(p, q)\in\mathfrak{F}$. Note that $T_qT_p<G(p, q)$ and hence $T_qT_p\in\mathfrak{F}$, $T_q\triangleleft T_qT_p$ and $C_{T_qT_p}(T_q)=T_q$. Thus $(q, p)\in E(\Gamma_s(\mathfrak{F}))$. It means that if a hereditary formation $\mathfrak{F}$ is recognized by $\Gamma_s$ then $\Gamma_s(\mathfrak{F})$ is undirected.

$(b)$ Let $G$ be a $s$-critical group for $\mathfrak{F}_\Gamma=\mathfrak{F}$. Since $\mathfrak{F}$ is hereditary, if $V(\Gamma(G))\not\subseteq V(\Gamma(\mathfrak{F}))$ then  $G\simeq Z_q$ for some $q\not\in\pi(\mathfrak{F})$.

 Assume now that $V(\Gamma(G))\subseteq V(\Gamma(\mathfrak{F}))$. So $E(\Gamma(G))\not\subseteq E(\Gamma(\mathfrak{F}))$. Note that $(p, q)\in\Gamma_s(G)$ if and only if $(p, q)\in\Gamma_s(N_G(P))$ for some Sylow $p$-subgroup $P$ of $G$.

  Since $G$ is $s$-critical for $\mathfrak{F}$, it follows that $G=N_G(P)$ for some Sylow subgroup $P$ of $G$ and there is $q\in\pi(G)$ such that $(p, q)\not\in\Gamma(\mathfrak{F})$. Hence $G/C_G(P)$ has a subgroup $A/C_G(P)\simeq Z_q$. Note that there is a $q$-element $x\in G$ such that $\langle xC_G(P)\rangle=A/C_G(P)$.
 Now $(p, q)\in E(P\langle x\rangle)$. Therefore $G=P\langle x\rangle$ and $\pi(G)=\{p, q\}$ .

Suppose that $G$ has a non-nilpotent proper subgroup. Then there is a Schmidt subgroup $K$ of $G$. Note that all subgroups of $G$ has normal Sylow $p$-subgroups. Therefore $K$ is a Schmidt group with the normal Sylow $p$-subgroup and $K\not\in \mathfrak{F}_\Gamma$, a contradiction. Thus   $G$ is a Schmidt group. Hence $\mathfrak{F}_\Gamma$ is a hereditary formation with the Shemetkov property.      According to corollary 2.4.3 \cite[p. 94]{s7} $\mathfrak{F}_\Gamma$ is a solubly saturated formation. $\square$

\smallskip

Recall \cite{g6} that a local formation $\mathfrak{F}=LF(F)$ where $F(p)=\mathfrak{S}_{\pi(F(p))}$ for all $p\in\pi(\mathfrak{F})$ and $F(p)=\emptyset$ otherwise is called a covering formation if $p\in\pi(F(p))$ and $p\in\pi(F(q))$ implies $q\in\pi(F(p))$. Such formations $\mathfrak{F}$ are hereditary saturated formations that contain every soluble group $G$ whose normalizers of Sylow subgroup are $\mathfrak{F}$-groups.

\begin{thm}\label{t4} Let $\mathfrak{X}$ be a class of groups such that $\Gamma_s(\mathfrak{X})$ is undirected. Then $\mathfrak{X}_{\Gamma_s}\cap\mathfrak{S}$ is a covering formation of soluble groups. In particular $\mathfrak{X}_{\Gamma_s}\cap\mathfrak{S}$ is a hereditary formation. \end{thm}

\emph{Proof.} Assume that $\mathfrak{X}_{\Gamma_s}\cap\mathfrak{S}$ is not hereditary. It means that there exists a group $G\in \mathfrak{X}_{\Gamma_s}\cap\mathfrak{S}$ with a subgroup $H$ such that $(p, q)\in E(\Gamma_s(H))$ and $(p, q)\not\in E(\mathfrak{X}_{\Gamma_s})$. Note that also $(q, p)\not\in E(\mathfrak{X}_{\Gamma_s})$.

Let $P$ be a Sylow $p$-subgroup of $H$. There is a Hall $\{p, q\}$-subgroup   $H_1$ of $H$ that contains a Hall $\{p, q\}$-subgroup of $N_H(P)$. Now $H_1$ is non-nilpotent.

Let $G_1$ be a Hall $\{p, q\}$-subgroup of $G$ that contains $H_1$. Since $(p, q)\not\in E(\Gamma_s(G))$   and  $(q, p)\not\in E(\Gamma_s(G))$, we see that $\Gamma_s(G_1)$ is two isolated vertices $p$ and $q$.

Assume that there exist non-nilpotent biprimary groups  $K$ such that  $\Gamma_s(K)$ is two isolated vertices $p$ and $q$. Let us chose  the minimal order group $K$ among them. Since $\mathfrak{N}$ is a saturated formation and $\Gamma_s$ is $\textbf{Q}$-closed we see that $K$ has an unique minimal normal subgroup $N$. Since $K$ is soluble, without lose of generality we may assume that $N$ is a $p$-group.   Now $K=N\leftthreetimes M$ where $M\in\mathfrak{N}$. Let $P$ a Sylow $p$-subgroup and $Q$ be a Sylow $q$-subgroup of $M$. Now $NP$ is a normal subgroup of $K$. Since $K$ is not nilpotent, $Q\not\leq C_K(NP)$. Thus $\Gamma_s(K)=(p, q)$, a contradiction.

Now $G_1$ is nilpotent and contains a non-nilpotent subgroup $H_1$, the contradiction. Thus  $\mathfrak{X}_{\Gamma_s}\cap\mathfrak{S}$ is  hereditary. So  $G\in \mathfrak{X}_{\Gamma_s}\cap\mathfrak{S}$ if and only if $G\in\mathfrak{S}$ and all its Schmidt subgroups are $\mathfrak{F}$-groups. It means that $(\mathfrak{X}_{\Gamma_s}\cap\mathfrak{S})_{\Gamma_{Sch}}\cap \mathfrak{S}=\mathfrak{X}_{\Gamma_s}\cap\mathfrak{S}$.  By theorem 4.4 $\mathfrak{X}_{\Gamma_s}\cap\mathfrak{S}$ is saturated. Thus $\mathfrak{X}_{\Gamma_s}\cap\mathfrak{S}$ is hereditary saturated  formation with the Schemetkov property in the soluble universe.

It is known \cite[p. 105]{s7} that in this case  $\mathfrak{X}_{\Gamma_s}\cap\mathfrak{S}$ has  a local definition $F$ such that $p\in\pi(F(p))$ and $F(p)=\mathfrak{S}_{\pi(F(p))}$ for all $p\in\pi( \mathfrak{X}_{\Gamma_s}\cap\mathfrak{S})$   and $F(p)=\emptyset$ otherwise.

Since $\Gamma_s(\mathfrak{X})$ is undirected, if a Schmidt $(p, q)$-group belongs to $\mathfrak{F}$  then   every Schmidt $(q, p)$-group  also belong to $\mathfrak{F}$. It means that if $q\in \pi(F(p))$ then $p\in\pi(F(q))$.     Thus  $\mathfrak{X}_{\Gamma_s}\cap\mathfrak{S}$ is a covering formation of soluble groups. $\square$


\begin{cor}\label{4.7} Let $\mathfrak{F}$ be a hereditary formation. If  $\mathfrak{F}$ is recognized by $\Gamma_s$    then the following statements holds:

$(a)$ $\Gamma_{s}(\mathfrak{F})$ is undirected.

$(b)$  $\mathfrak{F}$ is a solubly saturated formation with the Schemetkov property.

$(c)$ $\mathfrak{F}\cap\mathfrak{S}$ is a covering formation of soluble groups.  \end{cor}

\smallskip

\textbf{Problem 3.}  Assume that a hereditary formation $\mathfrak{F}$ satisfies (a)-(b) from corollary \ref{4.7}. Is $\mathfrak{F}$ recognized by $\Gamma_s$?

\section{Applications}

\subsection{Applications of the Hawkes graph}

Recall that Hawkes \cite{g5} show that if $\Gamma_H(G)$ has  no cycles then $G$ is a Sylow tower group for a linear ordering $\phi$   of  $\mathbb{P}$, i.e if $\{p_1,\dots, p_n\}$ are prime divisors of $G$ such that $p_i\leq_\phi p_j$   for $i<j$ then $G$ has normal Hall $\{p_1,\dots, p_i\}$-subgroups for all $i\leq n$.

\begin{thm} Let $\mathfrak{F}$ be a class of groups.  If $\pi\subseteq\pi(\mathfrak{F})$ and there are no edges from vertices from $\pi'$ to vertices from $\pi$ in $\Gamma_H(\mathfrak{F})$ then every group of $\mathfrak{F}$ has the normal Hall $\pi$-subgroup.
  \end{thm}


\emph{Proof.} 
Without lose of generality we may assume that $\mathfrak{F}=\mathfrak{F}_{\Gamma_H}$. Let a group $G$ be a counterexample of minimal order to  theorem. Let $N$ be a minimal normal subgroup of a group $G$.

Assume that $N$ is a $\pi$-group. Since $G/N\in\mathfrak{F}$, we see that $G/N$  has the normal Hall $\pi$-subgroup $H/N$. Now $H$ is  the normal Hall $\pi$-subgroup of $G$, a contradiction.

Assume now that $N$ is not a $\pi$-group. From the definition of $\Gamma_H(\mathfrak{F})$ and $C_G(N)\geq \mathrm{O}_{p', p}(G)$ for all $p\in\pi(N)$  it follows that $\pi(G/C_G(N))\cap\pi=\emptyset$. It means that $C_G(N)\triangleleft G$ contains all $\pi$-elements of $G$.

If $C_G(N)<G$ then there is a Hall $\pi$-subgroup $H\triangleleft C_G(N)\triangleleft G$. Hence $H$ is the normal Hall $\pi$-subgroup of $G$, a contradiction.

Suppose that $C_G(N)=G$. Then $N$ is an elementary abelian $p$-subgroup of $G$ for $p\not\in\pi$.    Since $G/N\in\mathfrak{F}$, we see that $G/N$  has the normal Hall $\pi$-subgroup $HN/N$ where $H$ is a $\pi$-group of $G$.  Now $N$ is the normal Hall $p$-subgroup of $HN\triangleleft G$. By Schur-Zassenhaus theorem $N$ has a complement $K$ in $HN$. Since $C_{HN}(N)=HN$, we see that $K$ is the normal Hall $\pi$-subgroup of $G$, the final  contradiction. $\square$

\begin{cor} Let $\mathfrak{F}$ be a class of groups, $V(\Gamma_H(\mathfrak{F}))=\pi_1\cup\pi_2$ where $\pi_1\cap\pi_2=\emptyset$ and there are no edges between $\pi_1$ and $\pi_2$. Then every group in $\mathfrak{F}$ is a direct product of its Hall $\pi_1$-subgroup and   $\pi_2$-subgroup. \end{cor}

\begin{thm} Let $A, B$ and $C$ be subgroups of a solvable group $G$ with pairwise coprime indexes. Then $\Gamma_H(G)=\Gamma_H(A)\cup\Gamma_H(B)\cup\Gamma_H(C)$.\end{thm}

\emph{Proof.} Let $\Gamma=\Gamma_H$. From  proposition \ref{p1} it follows that $\Gamma(A)\cup\Gamma(B)\cup\Gamma(C)\subseteq \Gamma(G)$. Assume that $ \Gamma(G)\not\subseteq\Gamma(A)\cup\Gamma(B)\cup\Gamma(C)$. Note that $V(\Gamma(A)\cup\Gamma(B)\cup\Gamma(C)) =V(\Gamma(G))$.
Now there exists $(p, q)\in E(\Gamma(G))$ such that  $(p, q)\not\in E(\Gamma(A)\cup\Gamma(B)\cup\Gamma(C))$.

Note that for every $p\in\pi(G)$ at least two subgroups from $A, B$ and $C$ contain a Sylow $p$-subgroup of $G$. It means that for every $\{p, q\}\subseteq\pi(G)$  at least one subgroup from $A, B$ and $C$ contains a Hall  $\{p, q\}$-subgroup of $G$. Since $G$ is soluble, it follows that  all Hall $\{p, q\}$-subgroups of $G$ are conjugate in $G$.

Assume that $p\neq q$. It means that one of subgroups $A, B$ and $C$ contains a Hall    $\{p, q\}$-subgroup $H$ of $G$.  Let $Q$ be a Sylow $q$-subgroup of $H$. By assumption $(p, q)\not\in\Gamma(H)$. Let us mention that $(p, q)\in E(\Gamma(H))$ if and only if $(p, q)\in E(\Gamma(H/\mathrm{O}_{p'}(H)))$. It means that $H/\mathrm{O}_{q}(H)$ is a $p$-group.   Now $Q$ is normal in $H$.     From $(p, q)\in\Gamma(G)$ it follows that
  \[1\neq Q\mathrm{O}_{p'}(G)/\mathrm{O}_{p'}(G)\triangleleft H\mathrm{O}_{p'}(G)/\mathrm{O}_{p'}(G).\]

  From $\mathrm{F}(G/\mathrm{O}_{p'}(G))=\mathrm{O}_{p}(G/\mathrm{O}_{p'}(G))$ it follows that \[C_{G/\mathrm{O}_{p'}(G)}(\mathrm{O}_p(G/\mathrm{O}_{p'}(G)))=C_{G/\mathrm{O}_{p'}(G)}(\mathrm{F}(G/\mathrm{O}_{p'}(G)))\subseteq \mathrm{F}(G/\mathrm{O}_{p'}(G))=\mathrm{O}_p(G/\mathrm{O}_{p'}(G)).\]

  Note that $\mathrm{O}_p(G/\mathrm{O}_{p'}(G))\leq H\mathrm{O}_{p'}(G)/\mathrm{O}_{p'}(G)$. Now \[Q\mathrm{O}_{p'}(G)/\mathrm{O}_{p'}(G)\subseteq C_{G/\mathrm{O}_{p'}(G))}(\mathrm{O}_{p}(G/\mathrm{O}_{p'}(G)))\subseteq \mathrm{O}_{p}(G/\mathrm{O}_{p'}(G)).\] So $Q\mathrm{O}_{p'}(G)/\mathrm{O}_{p'}(G)=1$, a contradiction.

Assume now that $q=p$.  Let $r\in\pi(G)$. From \[\mathrm{F}(G/\mathrm{O}_{p'}(G))=\mathrm{O}_p(G/\mathrm{O}_{p'}(G))  \,\,  and \,\,  C_{G/\mathrm{O}_{p'}(G)}(\mathrm{F}(G/\mathrm{O}_{p'}(G)))\subseteq\mathrm{F}(G/\mathrm{O}_{p'}(G))\] it follows that $\mathrm{O}_{r}(H\mathrm{O}_{p'}(G)/\mathrm{O}_{p'}(G))=1$ for every Hall    $\{p, r\}$-subgroup $H$ of $G$.

Note that one of subgroups $A, B$ and $C$ contains a Hall    $\{p, r\}$-subgroup $H$ of $G$.
So $(p, p)\not\in\Gamma(H)$. It means $(H\mathrm{O}_{p'}(G)/\mathrm{O}_{p'}(G))/\mathrm{O}_{p', p}(H/\mathrm{O}_{p'}(G))$ is a $r$-group. From $\mathrm{O}_r(H\mathrm{O}_{p'}(G)/\mathrm{O}_{p'}(G))=1$ it follows that $H/\mathrm{O}_{p'}(G)$ has the normal Sylow $p$-subgroup $P\mathrm{O}_{p'}(G)/\mathrm{O}_{p'}(G)$.

Thus for every  $r\in\pi(G)$ there exists a Hall  $\{p, r\}$-subgroup $H$ of $G$ such that $P\mathrm{O}_{p'}(G)/\mathrm{O}_{p'}(G)\triangleleft H\mathrm{O}_{p'}(G)/\mathrm{O}_{p'}(G)$. Now $P\mathrm{O}_{p'}(G)/\mathrm{O}_{p'}(G)\triangleleft G/\mathrm{O}_{p'}(G)$. Hence $(p, p)\not\in\Gamma(G)$, the contradiction. $\square$

Note that the condition of coprime indexes can't be omitted. Consider the symmetric group of degree 4. It is  the product of any two of the following subgroups: its Sylow $2$-subgroup, the alternating group of degree 4 and a subgroup that is isomorphic to the symmetric group of degree 3. The union of their Hawkes graphs is $\{ (2, 3), (3, 2)\}$. But the Hawkes graph of the symmetric group of degree 4 is $\{ (2, 3), (3, 2), (2, 2)\}$.

\begin{cor} Let $\mathfrak{F}$ be a formation such that $\mathfrak{F}$ is recognized by $\Gamma_H$. If a group $G$ contains three soluble $\mathfrak{F}$-subgroups $A$, $B$ and $C$ with pairwise coprime indexes      then $G\in\mathfrak{F}.$  \end{cor}

\emph{Proof.} Since $G$ contains three soluble  subgroups  with pairwise coprime indexes,  $G$ is soluble by Wielandt's theorem. Now $\Gamma_H(G)=\Gamma_H(A)\cup\Gamma_H(B)\cup\Gamma_H(C)=\Gamma_H(A\times B\times C)$. Since $A\times B\times C\in\mathfrak{F}$ and $\mathfrak{F}$ is recognized by $\Gamma_H$, we see that $G\in\mathfrak{F}.$              $\square$

\subsection{Applications of Schmidt graph}

Here we will describe Schmidt graphs of minimal simple non-abelian  groups.

\begin{thm}[\textbf{Thompson} {\cite[p. 190]{s10}}]\label{tT}  All minimal simple non-abelian groups are:

$(1)$ $PSL(2,2^p)$ where $p$ is a prime.

$(2)$ $PSL(2,3^p)$ where $p$ is an odd prime.

$(3)$ $PSL(2,p)$ where $p>5$ is a prime and $p^2+1\equiv 0 \mod 5$.

$(4)$ $Sz(2^p)$ where $p$ is an odd prime.

$(5)$ $PSL(3,3)$.
\end{thm}

\begin{thm}[\textbf{Dickson} {\cite[p. 213]{s10}}]\label{tD}   Any subgroup of $PSL(2, p^n)$ is isomorphic to a group from one of the following families of groups.

$(a)$ Elementary abelian $p$-groups.

$(b)$ Cyclic groups $Z_m$ of order $m$, where $m$ is a divisor of $(p^n\pm 1)/d$ and $d=(p-1,2)$.

$(c)$ Dihedral groups of order $2m$, where $m$ is defined in $(b)$.

$(d)$ Alternating group $A_4$ if $p>2$ or $p=2$ and $n\equiv 0 \mod 2$.

$(e)$ Symmetric group $S_4$ if $p^{2n}\equiv 1 \mod16$.

$(f)$  Alternating group $A_5$ if $p=5$ or  $p^{2n}\equiv 1 \mod 5$.

$(g)$  A semi-direct product of an elementary abelian $p$-group of order $p^m$ and a cyclic group of order $k$, where $k$ is a divisor of $p^m-1$ and $p^n-1$.

$(h)$  The group $PSL(2,p^m)$ if $m$ is a divisor of $n$, or the group $PGL(2,p^m)$ if $2m$ is a divisor of $n$.
\end{thm}

Let $H=N\leftthreetimes C$ be a subgroup of type (g). Then it is straightforward to check that $C_H(N)=N$. Let $G\simeq PSL(2, 2^p)$ be a minimal simple group. Note  that subgroups of types (f) and (h) are not proper subgroups of $G$.

 Direct calculations show that

 \begin{prop} The following statements hold:

  $(a)$ Let $p$ be a prime then $V(\Gamma_{Sch}(PSL(2, 2^p)))=\pi(2(2^{2p}-1))$ and \[E(\Gamma_{Sch}(PSL(2, 2^p)))=\{(2, q) | q\in\pi(2^p-1)\}\cup \{(q, 2) | q\in\pi(2^{2p}-1)\}.\]

  $(b)$ Let $p$ be an odd prime then $V(\Gamma_{Sch}(PSL(2, 3^p)))=\pi(3(3^{2p}-1))$ and \[E(\Gamma_{Sch}(PSL(2, 3^p)))=\{(3, q) | q\in\pi(3^p-1)\setminus\{2\}\}\cup \{(2, 3)\}\cup \{(q, 2) | q\in\pi(3^{2p}-1)\setminus\{2\}\}.\]

   $(c)$ Let $p$ be an odd prime such that $p^2+1\equiv 0 \mod 5$ then $V(\Gamma_{Sch}(PSL(2, p)))=\pi(p(p^{2}-1))$  \[E(\Gamma_{Sch}(PSL(2, p)))=\{(2, p) | q\in\pi(p-1)\setminus\{2\}\}\cup \{(2, 3)\}\cup \{(q, 2) | q\in\pi(p^2-1)\setminus\{2\}\}.\]

 \end{prop}

\begin{thm}[\textbf{Suzuki} {\cite{s11}}]\label{tSz} 
 Any subgroup of $G\simeq Sz(q)$ is isomorphic to a subgroup of one of the following  groups where $q=2^p$.

$(a)$ Frobenius groups of order $q^2(q-1)$, $H=QK$, $Q\triangleleft H$, $Q\in Syl_q(H)$ and $K$ is cyclic of order $q-1$.

$(b)$ Dihedral groups of order $q(q-1)$.

$(c)$ Cyclic groups $A_i$, $i=1,2$ of orders $q\pm r+1$  where $r^2=2q$.

$(d)$   $B_i=N_G(A_i)$ of order  $4(q\pm r+1)$.

$(e)$ $Sz(s)$ if $q$ is a power of $s$.
\end{thm}

 \begin{prop}

  Let $p$ be an odd prime then $V(\Gamma_{Sch}(PSL(3,3)))=\pi(2(2^{2p}+1)(2^p-1))$ and \[E(\Gamma_{Sch}(Sz(2^p)))=\{(2, q) | q\in\pi(2^p-1)\}\cup \{(q, 2) | q\in\pi((2^p-1)(2^{2p}+1))\}.\]
\end{prop}

 \begin{prop}

   For $PSL(3, 3)$ holds $V(\Gamma_{Sch}(Sz(2^p)))=\{2, 3, 13\}$ and \[E(\Gamma_{Sch}(PSL(3,3)))=\{(2, 3), (3, 2), (13, 3)\}.\]
\end{prop}

\begin{thm}\label{sol} Let $G$ be a group. The following statements holds.

$(a)$ If $\Gamma_{Sch}(G)$ does not contain cycles then $G$ is soluble.

$(b)$ Let $G$ be a group. If every cycle of $\Gamma_{Sch}(G)$ does not contain edge $(2, q)$ where $q\in\pi(2^p-1)$ for some prime $p$     then $G$ is soluble.

$(c)$  Let $G$ be a group. If every cycle of $\Gamma_{Sch}(G)$ has length greater then 3  then $G$ is soluble.
 \end{thm}

\emph{Proof.} (a) Let $\Gamma=\Gamma_{sch}$. Assume that $\Gamma(G)$ does not contain cycles and $G$ is not soluble. Then $G$ contains $s$-critical for $\mathfrak{S}$ subgroup $H$. It is well known that $H/\Phi(H)$ is minimal simple non-abelian group.

 So if $G$ is not soluble then $\Gamma(G)$ contains $\Gamma(H)$  for some minimal simple non-abelian group $H$.

Let $H\simeq PSL(2, 2^p)$   where $p$ is a prime. Then all cycles of $\Gamma(H)$ are $\{(2, q), (q, 2)\}$ where $q\in\pi(2^p-1)$.

Let $H\simeq PSL(2,3^p)$ where $p$ is an odd prime. Since $p$ is odd, $3^p-1\equiv 2 \mod 4$. Hence $|\pi(3^p-1)|>1$.   Now all cycles of $\Gamma(H)$ are $\{(2, 3), (3, q), (q, 2)\}$ where $q\in\pi(3^p-1)\setminus\{2\}$.

Let $H\simeq PSL(2,p)$ where $p>5$ is a prime and $p^2+1\equiv 0 \mod 5$. Now $3$ divides $p^2-1.$  So all cycles of $\Gamma(H)$ are $\{(2, 3), (3, 2)\}$ and $\{(2, p), (p, q), (q, 2)\}$ where $q\in\pi(p-1)\setminus\{2\}$.

Let $H\simeq Sz(2^p)$ where $p$ is an odd prime. Then all cycles of $\Gamma(H)$ are $\{(2, q), (q, 2)\}$ where $q\in\pi(2^p-1)$.

Let $H\simeq PSL(3,3)$. Then the only cycle of $\Gamma(H)$ is $\{(2, 3), (3, 2)\}$.

Thus $H$ and hence $G$ contains a cycle, the contradiction. So $G$ is soluble.

Let us prove $(b)$ and $(c)$.   As follows from the proof of $(a)$ every cycle of minimal simple group has an edge $(2, q)$ where $q\in\pi(2^p-1)$ for some prime $p$ and a length at most 3. $\square$


\begin{thm} Let $G$ be a group. If $\Gamma_{Sch}(G)$ has no cycles then $G$ has a Sylow tower for some linear ordering $\phi$ of $\mathbb{P}$.    \end{thm}

\emph{Proof.}  Assume that the theorem is wrong and let a group $G$ be a minimal order counterexample. From (a) of theorem \ref{sol} it follows that $G$ is soluble. Since $\Gamma_{Sch}(G)$ has no cycles, it follows that there is $p\in\pi(G)$ such that there no edges to $p$ in  $\Gamma_{Sch}(G)$. Consider a Sylow $p$-subgroup $P$ of $G$. Since $G$ is soluble, we see that for every $q\in\pi(G)$ there exists a Hall $\{p, q\}$-subgroup $Q$ of $G$ such that $P\leq Q$. Note that  $\Gamma_{Sch}(Q)$ is either two disjoint points or an arrow from $p$ to $q$. In the first case $Q$ is nilpotent and hence $Q\leq N_G(P)$. In the second case $Q$ is $q$-nilpotent and again $Q\leq N_G(P)$. Thus $P\triangleleft G$.

Note that $\Gamma_{Sch}(G/P)\subseteq\Gamma_{Sch}(G)$ and $|G/P|<|G|$. Hence $G/P$ has a Sylow tower for some linear ordering $\phi$ of $\pi(G/P)$. If we extend $\phi$ to $\pi(G)$ such that $p$ be the greatest element of $\pi(G)$ then $G$ will have the Sylow tower for this ordering, the contradiction. $\square$

\smallskip



\begin{cor}[Hawkes \cite{g5}] Let $G$ be a group. If $\Gamma_{H}(G)$ has no  cycles then $G$ has a Sylow tower for some linear ordering $\phi$ of $\mathbb{P}$.    \end{cor}

The following theorem we will prove with the help of some results that was obtained modulo the classification of  finite simple non-abelian groups.

\begin{thm}\label{sct} Let $\mathfrak{F}$ be a class of groups, $V(\Gamma_{Sch}(\mathfrak{F}))=\pi_1\cup\pi_2$ where $\pi_1\cap\pi_2=\emptyset$ and there are no edges between $\pi_1$ and $\pi_2$. Then every $\mathfrak{F}$-group   is a direct product of its Hall $\pi_1$-subgroup and   $\pi_2$-subgroup. \end{thm}

\emph{Proof.} We need to prove the theorem only for $\mathfrak{F}_{\Gamma_{Sch}}$. Suppose that the theorem is false and let a group $G$ be a counterexample of minimal order and $N$ be a minimal normal subgroup of $G$.

Assume that  $N$ is an abelian $p$-group. It is clear that $N\neq G$. Since $\Gamma_{Sch}$ is $\textbf{Q}$-closed, we see that $G/N=H_1N/N\times H_2N/N$ where $\pi(H_iN/N)\subseteq \pi_i$ for $i=1,2$. Without lose of generality we may assume that $p\in\pi_1$. Now $G$ has the normal Hall $\pi_1$-subgroup $H_1N$. Note that $N$ is the normal Hall $p$-subgroup of $H_2N$. By Schur-Zassenhaus theorem $N$ has a complement in $NH_2$. That is why we may assume that $H_2$ is a Hall $\pi_2$-subgroup of $H_2N$ and $G$. Note that $H_1N\triangleleft G$. Let $R$ be a Sylow $r$-subgroup of $H_2$ and $T$ be the set of all Sylow $q$-subgroups of $H_1N$. We see that  the order of $T$ is $\pi_1$-number. Since $NH_1\triangleleft G$, we see that $R$ acts on $T$ and size of any orbit is $r^\alpha$ where $\alpha\geq 0$. From $(r, q)=1$ it follows that $R$ has a fixed point in $T$. So $R\leq N_{G}(Q)$ for some Sylow $q$-subgroup $Q$ of $G$. Since $q\in\pi_1$ and $r\in\pi_2$ by the statement of theorem $R\leq C_G(Q)$. By analogy we can show that for every $q\in\pi_1$ there is a Sylow $q$-subgroup $Q$ of $NH_1$ such that $R\leq C_G(Q)$. Hence $R\leq C_G(H_1N)$. So $R\leq C_G(H_1N)$ for every Sylow subgroup $R$ of $H_2$. It means that $H_2\leq C_G(H_1N)$. Thus $G=H_2\times H_1N$, a contradiction.

Now assume that $N$ is non-abelian. In \cite{kaz} Kazarin and etc. proved that the Sylow graph of an almost simple group $H$ is connected, i.e. if $p, q\in\pi(H)$ then there are primes $p=p_0,\dots, p_n=q$ such that either $(p_{i-1}, p_i)\in E(\Gamma_s(H))$ or $(p_{i}, p_{i-1})\in E(\Gamma_s(H))$ for every $i=1,\dots, n$. From $\Gamma_s(G)\subseteq\Gamma_{Sch}(G)$ it follows that $G\neq N$. By our assumption $G/N=H_1N/N\times H_2N/N$. Assume that $\pi(N)\cap\pi_i\neq\emptyset$ for $i=1,2$. Since $\Gamma_{Sch}(N)$ is connected and $\Gamma_{Sch}(NH_i)\subseteq\Gamma_{Sch}(G)$ for $i=1,2$, we see that $\pi_1$ and $\pi_2$ are connected in $\Gamma_{Sch}(G)$, a contradiction. Thus $N$ is ether $\pi_1$-group or $\pi_2$-group. Now the proof is the same as in case when $N$ is abelian. $\square$


The following result follows from the first part of the   proof of theorem \ref{sct}.

\begin{cor}\label{sct1} Let $\mathfrak{F}$ be a class of soluble groups, $V(\Gamma_s(\mathfrak{F}))=\pi_1\cup\pi_2$ where $\pi_1\cap\pi_2=\emptyset$ and there are no edges between $\pi_1$ and $\pi_2$. Then every $\mathfrak{F}$-group   is a direct product of its Hall $\pi_1$-subgroup and   $\pi_2$-subgroup. \end{cor}

\end{document}